\documentclass[a4paper, 11pt]{amsart}

\usepackage{amssymb,amsmath,amsthm,amsfonts,amscd,mathrsfs}
\usepackage[all]{xy}
\usepackage[dvips]{graphicx}
\usepackage{hyperref} 

\def\Spec{\operatorname{Spec}}

\def\Hom{\operatorname{Hom}}
\def\dim{\operatorname{dim}} 
\def\codim{\operatorname{codim}} 
\def\trop{\operatorname{trop}} 
\def\red{\operatorname{red}} 

\def\pr{\operatorname{pr}} 
\def\Star{\operatorname{Star}} 

\def\ker{\operatorname{ker}} 

\def\val{\operatorname{val}} 
\def\charact{\operatorname{char}} 
\def\Pic{\operatorname{Pic}} 

\def\bbA{\mathbb A}

\def\bbK{\mathbb K}

\def\bbP{\mathbb P}
\def\bbQ{\mathbb Q}

\def\bbZ{\mathbb Z}

\def\cG{\mathcal G}

\def\cO{\mathcal O}

\def\cT{\mathcal T}

\def\sN{\mathscr N}

\def\fX{\mathfrak X}
\def\fY{\mathfrak Y}

\def\fm{\mathfrak m}

\def\oY{\overline{Y}}

\def\tY{\widetilde{Y}}
\def\tN{\widetilde{N}}
\def\tM{\widetilde{M}}
\def\tm{\widetilde{m}}
\def\Ssig{S_{\sigma}}

\def\Osig{O_{\sigma}}
\def\Vsig{V_{\sigma}}

\newtheoremstyle{style1}
{} 
{} 
{\itshape} 
{} 
{\bfseries} 
{.} 
{ } 
{} 

\theoremstyle{style1}

\newtheorem{Theorem}{Theorem}[section]
\newtheorem{Proposition}[Theorem]{Proposition}
\newtheorem{Lemma}[Theorem]{Lemma}
\newtheorem{Corollary}[Theorem]{Corollary}
\newtheorem{Conjecture}[Theorem]{Conjecture}

\newtheoremstyle{style2}{}{}{}{}{\itshape}{.}{ }{}
\theoremstyle{style2}

\newtheorem{Definition}[Theorem]{Definition}
\newtheorem{Remark}[Theorem]{Remark}
\newtheorem{Remarks}[Theorem]{Remarks}

\newtheorem{Question}[Theorem]{Question}

\begin{document}

\title{On Tropical Compactifications}
\author{Mark Luxton and Zhenhua Qu}

\begin{abstract}
We give an affirmative answer to a conjecture proposed by Tevelev \cite{Tev07} in characteristic 0 case: any variety contains a sch\"on very affine open subvariety. Also we show that any fan supported on the tropicalization of a sch\"on very affine variety produces a sch\"on compactification. Using toric schemes over a discrete valuation ring, we extend tropical compatifications to the non-constant coefficient case.
\end{abstract}

\maketitle

\section{Introduction and Statement of Results}
Compactification problems play an important role in algebraic geometry. Tropical compactifications are introduced and developed by Tevelev in \cite{Tev07} as certain nice compactifications for subvarieties of tori. This is applied to the study of compactification of moduli spaces of del Pezzo surfaces in a subsequent paper \cite{HKT07}. In an effort of applying tropical compactification to prove the following conjecture \ref{mainconj} proposed in \cite{KT06}, we obtained some new results on tropical compactifications and are able to prove the $X(3,6)$ case as well. This paper collects these results and their generalizations to non-constant coefficient case as we believe they may have some general interests both in algebraic geometry and in tropical geometry.

\begin{Conjecture}[\cite{KT06}]\label{mainconj}
The Chow quotient compactification $\overline{X}(3,n)$ has log canonical singularity and is the log canonical model of $X(3,n)$ when $n=6,7,8$, where $X(r,n)$ is the moduli space of ordered $n$ hyperplanes in $\bbP^{r-1}$ in general position.
\end{Conjecture}

We recall some results from \cite{Tev07}. Let $Y\subset T$ be a subvariety of a torus, $X(\Delta)$ a toric variety containing $T$ and let $\oY$ be the closure of $Y$ in $X(\Delta)$, we say $\oY$ is a \emph{tropical compactification} if $\oY$ is proper and the structure map $T\times\oY\rightarrow X(\Delta)$ is flat and surjective. It is a remarkable observation by Tevelev that to obtain a tropical compactification, $\Delta$ is supported on $\trop(Y)$, the tropicalization of $Y$, which explains the name ``tropical compactification''. 

From birational geometry point of view, one often seeks a ``minimal" compactification, as in conjecture \ref{mainconj} we are looking for the log canonical compactification. Thus a natural question is:

\begin{Question}\label{mainques}
Can we obtain a log canonical compactification through a tropical compactification for a subvariety of a torus?
\end{Question}

We know quite a few about the answer, but not completely though. A log canonical compactification has two requirements, a singularity requirement namely that the pair $(\oY, B)$ has log canonical singularity where $B=\oY\backslash Y$ is the boundary divisor, and a minimality requirement namely $K_{\oY}+B$ is ample. Sch\"on compactifications are introduced to fulfill the singularity requirement \cite{Tev07}. We say that $\oY$ is a \emph{sch\"on compactification}, if it is tropical and the structure map is smooth. Thus sch\"on compactifications have toroidal singularities which is a little stronger than log canonical, but in this situation it is reasonable and easy to manage. 

Sch\"on varieties (i.e. varieties which admit a sch\"on compactification) have many nice properties. If $Y$ is sch\"on, then either $Y$ is log minimal or $Y$ is preserved by a non-trivial subtorus (\cite{HKT07}, 3.1). If $\oY$ is a sch\"on compactification such that for any toric oribt $O\subset X(\Delta)$, $\oY\cap O$ is connected if $\dim (\oY\cap O)>0$, then the link of 0 of the $\trop(Y)$ only has top  reduced rational homology (\cite{Hac08}). We prove two other properties for sch\"on varieties (see theorem \ref{main1} and \ref{main2}).  

When $Y\subset T$ is sch\"on, we have a complete answer to question \ref{mainques}. We say $Y$ is \emph{h\"ubsch} if $Y$ is sch\"on and a sch\"on compactification is the log canonical compactification. We have a necessary and sufficient condition for a sch\"on variety to be h\"ubsch in terms of the combinatoric of $\trop(Y)$.

\begin{Theorem}\label{critofhub}
Assume $Y$ is a sch\"on compactification, then $Y$ is h\"ubsch if and only if $\trop(Y)$ has a minimal fan structure $\Delta$, and for any cone $\sigma\in\Delta$, $\Star(\sigma)$ is not preserved by any translation.
\end{Theorem}

There are many examples of h\"ubsch varieties, a hyperplane complement with connected matroid structure, a hypersurface nondegenerate with respect to the Newton polytope, a generic complete intersection, $M_{0,n}$, $Y^n$ (moduli space of del Pezzo surfaces of degree $9-n$ for $n\leq 6$ or $n=7$, $\charact k\neq 2$, \cite{HKT07}) as well as $X(3,6)$ in conjecture \ref{mainconj}. Our main results are the following.

\begin{Theorem}\label{main1}
Let $k$ be an algebraic closed field of characteristic 0, then any variety over $k$ contains a sch\"on very affine open subvariety.
\end{Theorem}

This answers a conjecture of Tevelev proposed in \cite{Tev07}, in characteristic 0 case. Note that a (dominant) morphism of very affine varieties $X\rightarrow Y$ has a corresponding (surjective) morphism of their tropicalization $\trop(X)\rightarrow\trop(Y)$, the above theorem suggests that any tropical variety is dominated by a ``good'' one. However a more precise statement is needed in tropical geometry.

\begin{Theorem}\label{main2}
If $Y\subset T$ admits a sch\"on compactification, then any fan $\Delta$ supported on $\trop(Y)$ is tropical, and hence sch\"on.
\end{Theorem}

In general not every fan is tropical (\cite{ST08}, 3.10), thus the above theorem is a special property of sch\"on varieties. The rest of the paper is devoted to generalizing tropical compactifications to non-constant coefficient case which is probably more interesting in tropical geometry. In non-constant case we consider a very affine variety $Y$ defined over the field of Puiseux series over $k$, $\trop(Y)$ is then a polyhedral complex which may not admit a fan structure. However $\trop(Y)$ can be realized as the fibre of a map of sets $\cT(Y)\rightarrow\bbQ_{\geq 0}$ at 1, where $\cT(Y)$ is the closure of cone over $\trop(Y)$ in $N_{\bbQ}\oplus\bbQ$. $\cT(Y)$ does admit fan structures and we take it as a replacement of $\trop(Y)$ in the non-constant coefficient case. This leads us to consider relative toric varieties, or toric schemes over a discrete valuation ring. Using $\cT(Y)$ and toric schemes, we are able to generalize tropical, sch\"on and h\"ubsch compactifications to non-constant coefficient case. With some appropriate modification, theorem \ref{main1} and \ref{main2} remain true in this case.

In the final part of the introduction, we claim the following result and give only an outline of the proof due to the length of its full details. A careful treatment can be found in the first author's thesis \cite{Lux08}.

\begin{Theorem}
$X(3,6)$ is h\"ubsch.
\end{Theorem}

Recall that $X(r,n)$ is the moduli space of ordered $n$ hyperplanes in $\bbP^{r-1}$ in general position. Galfand-MacPherson correspondence identifies $X(r,n)$ with $G^o(r,n)/H$ where $H$ is the maximal torus of $\bbP^{r-1}$. The tropicalization of $X(3,6)$ is the tropical Grassmannian $\cG_{3,6}''$ studied in \cite{SS04} where an explicit description of the minimal fan structure on $\cG_{3,6}''$ is obtained. We recovered the same result using geometric tropicalization (see subsection 2.3). The key point is the study of the map $Y^6\rightarrow X(3,6)$ and the corresponding map of tropicalizations. A nice smooth compactification of $Y^6$ with s.n.c. boundary is obtained by Naruki using $D_4$ cross-ratios of $E_6$ \cite{Nar82}, see also \cite{Sek94, Sek00} for generalizations. The geometric tropicalization is a minimal fan and $Y^6$ is h\"ubsch (\cite{HKT07}). It turns out that $\oY^6$ contains $X(3,6)$ as an open subset and this geometric tropicalization refines the minimal fan of $\cG_{3,6}''$. Using this refined fan and $\oY^6$, we show it has smooth structure map, then by Theorem \ref{main2}, it also has a smooth structure map using the minimal fan on $\cG_{3,6}''$. Finally we verify that the minimal fan on $\cG_{3,6}''$ satisfies the assumption in theorem \ref{critofhub} and complete the proof.

This paper is organized in the following way. Some preliminaries on tropicalization and tropical compactification are included in section 2. We prove theorem \ref{main1} in section 3 and theorem \ref{main2} in section 4. Section 5 introduces toric schemes over a discrete valuation ring in order to extend our results to the non-constant coefficient case, which is done in section 6.

\bigbreak
\emph{Acknowlegements.} We would like to take the opportunity here to thank all whose conversation and advice are helpful to us , especially our thesis advisor S\'ean Keel for his continued support and advice. We also thank E. Katz and J. Tevelev for many helpful discussions.

\section{Preliminaries}
\subsection{Tropicalization} We recall the basics of tropical geometry, see for example \cite{EKL06, Spe05}. Let $k$ be an algebraically closed field of any characteristic, $\bbK$ be the algebraic closure of $k(t)$. We have a natural valuation $v:\bbK^*\rightarrow \bbQ$. 

\begin{Definition}
For a subscheme $Y\subset T\cong (\bbK^*)^n$, the \emph{tropicalization} of $Y$ is defined to be 
$$\trop(Y)=\{(v(y_1),\ldots,v(y_n))\in\bbQ^n|(y_1,\ldots,y_n)\in Y(\bbK)\}\subset\bbQ^n.$$
\end{Definition} 
It has several equivalent descriptions.

\begin{Theorem}[\cite{Spe05} 2.1.2]\label{tropset}
The following sets are equal:
\begin{enumerate}
\item $\trop(Y)$ defined as above.
\item The set $\left\{(u(t_1),\ldots,u(t_n))|u:\cO(Y)\rightarrow\bbQ\cup\{\infty\}\right\}$, where $u$ runs over all ring valuations $\cO(Y)\rightarrow\bbQ\cup\{\infty\}$ extending $v$.
\item The set of $w\in\bbQ^n$ such that $\mathrm{in}_w f$ is not a monomial for any $f\in I\backslash\{0\}$, or equivalently $\mathrm{in}_w Y\neq\emptyset$.
\end{enumerate}
\end{Theorem}

(2) in theorem \ref{defoftrop} is also called a BG-set (short for Bieri-Groves set, \cite{BG84}). $\trop(Y)$ is a polyhedral complex of pure dimension ($=\dim Y$). If $I$ is generated by a set of functions in $k[t_1^{\pm},\ldots,t_n^{\pm}]$, then we may think $Y$ is a subscheme of $T_k$, and we call this constant coefficient case. When we are in constant coefficient case, $\trop(Y)$ is conical, and admits a fan structure. In this section we consider the constant coefficient case only.

\subsection{Tropical Compactifications}

Now assume $Y$ is a closed (irreducible) subvariety of $T_k$, let $M$ be the lattice of characters of $T$, $N$ the dual of $M$, $\trop(Y)$ can be naturally considered in $N_{\bbQ}=N\otimes_{\bbZ}\bbQ$. Let $\Delta\subset N_{\bbQ}$ be a fan, denote $X(\Delta)$ the corresponding toric variety, and $\oY$ (or $\oY(\Delta)$ when one wants to specify the fan $\Delta$) the closure of $Y$ in $X(\Delta)$.

\begin{Definition}[\cite{Tev07}]\label{defoftrop} We say $\Delta$ is a \emph{tropical fan} or $\oY$ is a \emph{tropical compactification} if $\oY$ is compact and the structure map $T\times_k\oY\rightarrow X(\Delta)$ is flat and surjective. We say $\Delta$ is \emph{sch\"on fan} or $\oY$ is a \emph{sch\"on compactification} if it is tropical and furthermore the structure map is smooth.
\end{Definition}

The basic result about tropical compactifications is the following.

\begin{Theorem}[\cite{Tev07}]
Assume $Y\subset T$ is rigid (i.e. the stablizer of $Y$ in $T$ is trivial), then tropical compactification always exists. If $\Delta$ is a tropical fan, then $\Delta$ is supported on $\trop(Y)$ and any refinement of $\Delta$ is also tropical. If $Y\subset T$ admits a sch\"on compactification, then any tropical fan is sch\"on.
\end{Theorem}

\begin{Remark}\label{importantremark}
It is often useful to consider an arbitrary $Y$, not necessarily irreducible. The reason is that if $\oY$ is a tropical (or sch\"on) compactification of an irreducible variety $Y$, and $W\subset X(\Delta)$ is a toric orbit closure, then $\oY\cap W$ has flat (or smooth) surjective structure map in this toric variety, being the pullback of the structure map of $\oY$, but $\oY\cap O$ is usually not irreducible where $O\subset W$ is the open orbit (see also lemma \ref{lem2}). We can define tropical and sch\"on compactification of an arbitrary closed scheme $Y\subset T$ in the same way as in definition \ref{defoftrop}. Note that it remains true in this general case that $|\Delta|=\trop(Y)$ when $\Delta$ is tropical, however some results are not true anymore, for example it is not clear if tropical compactification exists when $Y$ is non-reduced. If $Y$ is a closed subscheme and $\oY$ is a sch\"on compactification, then $Y$ is necessarily reducible. Since $T\times\oY\rightarrow X(\Delta)$ is smooth, we conclude that $\oY$ is a disjoint union of irreducible components. Let $\oY'\subset \oY$ be an irreducible component, then there is a possibly smaller toric open subset $X'\subset X(\Delta)$ containing $\oY'$ with smooth surjective structure map. It follows that the fan of $X'$, being a subfan of $\Delta$, is supported on $\trop(\oY'\cap T)$.
\end{Remark}

\subsection{Geometric Tropicalization}

Geometric tropicalization computes the set $\trop(Y)$ in terms of a nice compactification of $Y$, namely a compactification with s.n.c. boundary divisors. Let $\oY$ be any such compactification, for each boundary divisor $D$, the valuation $\val_D$ restricted on $M$ determines a point in $N$, still denoted by $\val_D$. For each collection $S$ of boundary divisors with nonempty intersection, let $\sigma_S$ be the cone in $N$ generated by $\val_D$ with $D\in S$, then we have the following theorem.

\begin{Theorem}[\cite{HKT07}]\label{geomtrop}
$\trop(Y)$ is the union of all $\sigma_S$ where $S$ runs over all collection of boundary divisors with nonempty intersection.
\end{Theorem}

\begin{Remark}
The collection of $\sigma_S$ do not form a fan in genral. In the proof of theorem \ref{main1}, we search a nice compactification whose geometric tropicalization produces a fan with simplicial cones.
\end{Remark}

\begin{Remarks}\label{geomtroptoroidal}
In theorem \ref{geomtrop}, the requirement for the boundary being s.n.c. can be weekened. In fact, when $Y\subset\oY$ is a toroidal embedding without self-intersection (as defined in \cite{KKMSD}), theorem \ref{geomtrop} remains true. 

Mumford associates to every toroidal embedding $Y\subset\oY$ a conical polyhedral complex with intergral structure $\Delta=(|\Delta|,\sigma^S,M^S)$, whose cells are in one-to-one correspondence with the strata of $\oY$. Recall that for a strata $S$, $M^S$ is the group of Cartier divisors of $\Star(S)$ supported on $\Star(S)\backslash Y$, $M^S_+$ is the subgroups of $M^S$ consisted of effective Cartier divisors, $N^S=\Hom(M^S,\bbZ)$, and the cell $\sigma^S$ of $\Delta$ is spanned by $v\in N^S$ such that $\left\langle v,D\right\rangle\geq 0$ for any $D\in M^S_+$.

Assume $T=\Spec k[M]$, we have a morphism $M\rightarrow M^S$ for any stratum $S$, namely $m\mapsto (m)$ on $\Star(S)$, hence a corresponding dual map $N^S\rightarrow N$ and $|\Delta|\rightarrow N_{\bbQ}$. The image of $\sigma^S\subset N^S_{\bbQ}$ in $N_{\bbQ}$ is exactly $\sigma_S$, thus the image of $|\Delta|$ is $\bigcup_S\sigma_S$ as in Theorem \ref{geomtrop} (here we make no difference between a stratum $S$ and the collection of divisors containing $S$).

Let $\Delta'$ be a subdivision of $\Delta$ such that all cells are strictly simplicial. By the theory of toroidal embeddings, a subdivision corresponds to another toroidal embedding $Y\subset\oY'$ with a canonical birational morphism $\oY'\rightarrow\oY$. $\oY'$ has s.n.c. boundary divisors. By Theorem \ref{geomtrop}, image of $|\Delta'|$ is $\trop(Y)$, so is the image of $|\Delta|$.
\end{Remarks}

\subsection{Intrinsic Torus and H\"ubsch Varieties} 
For an affine variety $Y$, we say $Y$ is \emph{very affine} if it admits a closed embedding into some torus $T$. For different embeddings, the tropical varieties are different, but there is an embedding which has the richest tropical structure and dominates all others, namely the intrinsic one. Recall that for any variety $Y$, $\cO^*(Y)/k^*$ is a lattice. A closed embedding $Y\rightarrow T=\Spec k[M]$ is equivalent to a homomorphism $M\rightarrow\cO^*(Y)$ whose image generates $\cO(Y)$ as $k$-algebra. The intrinsic embedding uses exactly all the units, more precisely $Y\rightarrow T_{\textrm{in}}=\Spec k[M]$ is an intrinsic embedding if $M\rightarrow \cO^*(Y)/k^*$ is an isomorphism. Clearly any embedding $Y\rightarrow T$ factors through $Y\rightarrow T_{\textrm{in}}\rightarrow T$ and the tropicalization in $T$ is the image of the tropicalization in $T_{\textrm{in}}$ under the corresponding map $N_{T_{\textrm{in}}}\rightarrow N_T$. 

\begin{Definition}
Given $Y\subset T$, we say $Y$ is \emph{sch\"on} in $T$ if it admits a sch\"on compactification. We say a very affine variety $Y$ is \emph{sch\"on} if it is sch\"on in the intrinsic torus.
\end{Definition}

\begin{Definition}
Given $Y\subset T$, we say $Y$ is \emph{h\"ubsch} in $T$ if it is sch\"on in $T$ and a sch\"on compactification is the log canonical compactification. We say a very affine variety $Y$ is \emph{h\"ubsch} if it is h\"ubsch in the intrinsic torus.
\end{Definition}

$Y$ being h\"ubsch in a torus $T$ imposes strong restriction on the combinatorics of its tropicalization.

\begin{Theorem}[\cite{HKT07}]
If $Y\subset T$ is h\"ubsch, then $\trop(Y)$ has a minimal fan structure $\Delta$, and $\oY(\Delta)$ is the log canonical compactification.
\end{Theorem}

\begin{Remark}
In general, $\trop(Y)$ does not have a minimal fan structure, see \cite{ST08} 5.2 for example.
\end{Remark}

\begin{Lemma}
If $Y$ is sch\"on (resp. h\"ubsch) in some torus $T$, then $Y$ is sch\"on (resp. h\"ubsch) in the intrinsic torus.
\end{Lemma}
\begin{proof}
Without loss of generality, we can assume $T_{\textrm{in}}\rightarrow T$ is surjective. Let $\Delta$ be a sch\"on fan supported on $\trop(Y)\subset N_T\otimes_{\bbZ}\bbQ$, and $\Delta'$ is any fan supportded on the intrinsic tropicalization $\trop(Y)\subset N_{T_{\textrm{in}}}\otimes_{\bbZ}\bbQ$ such that the induced map $\Delta\rightarrow\Delta'$ is a fan map. Let $\oY$ and $\oY'$ be closures in $X(\Delta)$ and $X(\Delta')$ respectively, we have a fibre diagram
$$ 
\begin{xymatrix}
{T_{\textrm{in}}\times_k \oY' \ar[r] \ar[d] & X(\Delta') \ar[d] \\
T\times_k \oY \ar[r] & X(\Delta). }
\end{xymatrix}
$$

Now suppose $Y\subset T$ is h\"ubsch. Let $\Delta$ be the canonical (minimal) fan supported on $\trop(Y)\subset N_{\bbQ}$. The compactification of $Y$ in $X(\Delta)$ is the log canonical model for $Y$, denoted by $\oY_{\textrm{lc}}$. Let $\Delta'$ be any fan supported on the intrinsic tropicalization, $\oY'$ be the closure of $Y$ in $X(\Delta')$, thus we have a canonical birational morphism $\oY'\rightarrow\oY_{\textrm{lc}}$. Recall that a ray $\rho$ in $\Delta'$ corresponds to some boundary divisors in $\oY'$ (namely $\oY'\cap V_{\rho}$) and a collection of rays form a fan if the corresponding divisors intersect, similarly for $\oY_{\textrm{lc}}$, thus the map $\oY'\rightarrow\oY_{\textrm{lc}}$ induces a fan map $\Delta'\rightarrow \Delta$. We may use $\oY_{\textrm{lc}}$ to get a fan structure on the intrinsic $\trop(Y)$ (remark \ref{geomtroptoroidal}), note that in Theorem \ref{geomtrop}, these cones do not necessarily form a fan, but here they do if we ignore repeated cones because the image of these cones in $N_{T}\otimes_{\bbZ}\bbQ$ form a fan, and 1-dimensional cones map to 1-dimensional cones. Call this fan $\Delta_{\textrm{lc}}$, then the closure of $Y$ in $X(\Delta_{\textrm{lc}})$ recovers $\oY_{\textrm{lc}}$.
\end{proof}

We state a criterion for $Y$ being h\"ubsch to close this section.

\begin{Theorem}[\cite{HKT07}, 9.1]
Assume $\oY\subset X(\Delta)$ is a sch\"on compactification, then it is the log canonical compactification iff for each toric orbit $\Osig\subset X(\Delta)$, $\oY\cap\Osig$ is log minimal (or equivalently, not preserved by a subtorus)
\end{Theorem}

\begin{Remark}
Theorem \ref{critofhub} follows from the above theorem and lemma \ref{lem1}.
\end{Remark}

\section{Embeddings into Toric Varieties}

We prove Theorem \ref{main1} in this section. Now assume $\charact k=0$. It suffices to prove the theorem for a smooth variety $Y$. Our strategy is to first compactify $Y$ so that the compactification $\oY$ is a smooth projective variety and the boundary $\oY\backslash Y$ is a s.n.c. divisor. We then study how to embed $\oY$ into a toric variety $X$ with one-to-one correspondence between strata of $\oY$ and toric strata of $X$ and with a smooth structure map. Certain requirements have to be imposed on the log structure of $Y$ (proposition \ref{schoncondition}). These requirements are achieved when we add more generic hyperplane sections to $Y$. 

The problem of embedding an arbitrary variety (possibly singular and nonprojective) into some toric variety is studied in \cite{Wlo93}. Embedding a variety into a given toric variety is equivalent to giving a compatible log structure, namely the pullback of the natural log structure on the toric variety. Our proof has a similar flavor. We also notice that in \cite{HKT07}, the authors obtained some similar requiements using quotient of affine conoid technique, but only applicable when the Picard group $\Pic \oY$ is a free abelian group of finite rank.

Let $\oY$ be a smooth proper variety over an algebraically closed field $k$ of characteristic 0. Let $D=\bigcup_{i\in I}D_i$ be a s.n.c. divisor. We have a stratification induced by $D$, a stratum is defined by $\bigcap_{j\in J}D_j-\bigcup_{i\in I\backslash J}D_i$ for a subset $J\subset I$, when $J=\emptyset$, the corresponding stratum is the open complement $Y=\oY\backslash D$. We don't require each stratum to be irreducible, this is a little different from the definition in \cite{KKMSD} where a stratum is an irreducible component of our stratum here, but it is more convenient in our purpose for avoiding repeated cones in the geometric tropicalization. For a stratum $S$, let $D_S$ be the set $\{D_i|S\subset D_i\}$, and let $\Star(S)$ be $Y-\bigcup_{D_i\notin D_S}$. Since $D$ is s.n.c., a stratum $S$ is regular (not necessarity irreducible), and also $\#D_S=\codim S$.

Fix $(M,\varphi)$, a pair of an abstract lattice $M$ and a group homomorphism $\varphi: M\rightarrow\cO^*(Y)$ such that $k[M]\rightarrow \cO(Y)$ is a surjection. For any submonoid $M'\subset M$, we have induced $k$-algebra homomorphism $\varphi_*: k[M']\rightarrow \cO(Y)$. For any stratum $S$, let $M_S\subset M$ be the submonoid of $M$ consisting all $m\in M$ such that $m$ is regular on $\Star(S)$ or equivalently $\val_{D_i} m\geq 0$ for any $D_i\in D_S$. By abuse of notation sometimes we use $m$ for which we actually mean $\varphi(m)$, for example $\val_D m$. Let $N$ be the dual of $M$, for any stratum $S$, let $\sigma_S\subset N_{\bbQ}$ be the dual cone of $M_S$.

\begin{Proposition}\label{schoncondition}
Notations as above, assume the following conditions are satisfied:
\begin{enumerate}
\item for any stratum $S$, $\Star(S)$ is affine and the induced map $k[M_S]\rightarrow\cO(\Star(S))$ is surjective, in particular, $Y$ is very affine.
\item for any stratum $S\neq Y$, and any $D_0\in D_S$, there exists $m\in M$ such that $\val_{D_0}m=1$ and $\val_{D_i}m=0$ for any $D_i\in D_S-{D_0}$.
\item the collection of cones $\{\sigma_S\}$ as $S$ runs over all strata is a fan $\Delta$.
\end{enumerate}
Then there is a canonical closed immersion $Y\rightarrow X(\Delta)$ and $S=Y\cap O_{\sigma_S}$ as scheme-theoretic intersection for any stratum $S$, hence the structure map is smooth and surjective, $X(\Delta)$ is a sch\"on compactification.
\end{Proposition}

\begin{proof}
First note that $\sigma_S$ is spanned by $\val_{D_i}$ for $D_i\in D_S$ since 
$$M_S=\left\{m\in M|\val_{D_i}m\geq 0\; \forall D_i\in D_S\right\}.$$ 
Condition (2) implies that $\val_{D_i}$ is part of a basis of $N$, i.e. $\sigma_S$ is strictly simplicial. 

For any stratum $S$, the surjection $k[M_S]\rightarrow\cO(\Star(S))$ and the fact that $\Star(S)$ is affine determines a closed embedding $\Star(S)\rightarrow X_{\sigma_S}$. We put an partial ordering on the strata, we say $S_1\leq S_2$ if $\Star(S_1)\subseteq\Star(S_2)$ or equivalently $D_{S_1}\subseteq D_{S_2}$ or equivalently $\sigma_{S_1}\leq\sigma_{S_2}$ ($\sigma_{S_1}$ is a face of $\sigma_{S_2}$). Thus the correspondence $S\mapsto\sigma_S$ is an isomorphism of partially ordered sets. For any strata $S_1,S_2$, let $S$ be the strata corresponding to $\sigma_{S_1}\cap\sigma_{S_2}$, we have a following commutative diagram,
$$\begin{xymatrix}
{\Star(S_1) \ar[d] & \Star(S) \ar[l] \ar[r] \ar[d] & \Star(S_2) \ar[d] \\
X_{\sigma_{S_1}} & X_{\sigma_S} \ar[l] \ar[r] & X_{\sigma_{S_2}} }
\end{xymatrix}$$
By gluing schemes and morphisms, we get a canonical closed embedding $Y\rightarrow X(\Delta)$.

Next we prove $S=Y\cap O_{\sigma_S}$ as scheme-theoretic intersection. Assume $D_S=\{D_1,\ldots,D_l\}$. Let $E_i=D_i\cap\Star(S)$. By condition (2), we can find $m_1,\ldots,m_l$ such that $\val_{D_i}m_j=\delta_{ij}$. Note that $m_i$ is part of a basis of $M$, the closed orbit $O_{\sigma_S}\subset X_{\sigma_S}=\Spec k[M_S]$ is defined by ideal $(m_1,\ldots,m_l)\subset k[M_S]$. Clearly $E_i=\Star(S)\cap (m_i=0)$ since $(m_i)=E_i$, we have $S=\bigcap E_i=\Star(S)\cap (m_1=\ldots=m_l=0)=\Star(S)\cap O_{\sigma_S}$.

The structure map $T_M\times Y\rightarrow X(\Delta)$ has fiber $S\times \ker(T_M\rightarrow O_{\sigma_S})$ over any point $x\in O_{\sigma_S}$, so it has reduced, equidimensional, regular fiber. It follows that the structure map is surjective and smooth.

\end{proof}

\begin{proof}[Proof of Theorem \ref{main1}]
Let $Y$ be a smooth variety, $\oY$ any smooth projective compactification of $Y$ with s.n.c. boundary divisor $D=\{D_1,\ldots,D_r\}$. Let $L$ be a very ample line bundle such that $L+D_i$ is very ample for any $i$. Choose a finite set $E_i\subset |L+D_i|$ of generic sections with $\#E_i\geq \dim |L+D_i|+\dim Y +1$ for $i=0,1,\ldots,r$ where we take $D_0$ to be the zero divisor. Let $E$ be the union of all $D_i$ and all divisors in $E_i$, $i=0,1,\ldots,r$. By Bertini's theorem, $E$ is s.n.c. Let $Y^o$ be the complement $\oY\backslash E$ and $M\cong\cO(Y^o)/k^*$. We show that $(\oY,E,M)$ satisfies the assumptions in proposition \ref{schoncondition}.

The following facts are frequently used: if $L$ is a very ample line bundle on a projective variety $Y$, and $s_1,\ldots,s_l$ are sections of $|L|$ in linear general position with $l=\dim |L|+1$, then the complement is very affine. Intersection of two very affine open subvarieties of $Y$ is again very affine, hence in the first statement, it is true for $l\geq\dim |L|+1$.

Let $S$ be any stratum, we have $\# D_S\leq \dim Y$. For any $D_i\in D\backslash D_S$, pair with a divisor $D_i'\in E_0\backslash D_S$. We can write 
$$E\backslash D_S=\bigcup_{i=0}^r F_i,$$
where $F_i=E_i\backslash D_S$ or $E_i\backslash D_S\cap\{D_i+D_i'\}$ if $D_i\in D\backslash D_S$. Each $F_i\subset |L+D_i|$ are sections in linear general position with $\#F_i\geq \dim|L+D_i|+1$, thus $\Star(S)=\cap(\oY\backslash F_i)$ is very affine. Since $\cO(\Star(S))$ is generated by all units, it is also generated by $M_S$ which contains all the units on $\Star(S)$.

To see condition (2), let $F\in D_S$ be any divisor, if $F$ is one of $D_i$'s, choose some $F'\in F_0\backslash D_S$ and consider $F+F'\in F_i$, otherwise $F\in F_i$ for some $i$, we choose another divisor $G\in F_i\backslash D_S$ and $G\neq F$ or $F+F'$. There is a unit $m\in M$ with the associated divisor $(m)= F-G$ or $F+F'-G$. In either case, this $m$ satisfies condition (2).

To verify condition (3), we show that for any two strata $S,S'$, their is a unit $m\in M$ such that $\val_F m>0$ for all $F\in D_S\backslash D_{S'}$ and $\val_F m\leq 0$ for all $F \in D_{S'}$. Indeed for each $F\in D_S\backslash D_{S'}$, as in the above argument we can find $m_F\in M$ such that $(m_F)=F-G$ (or $F+F'-G$ if $F\in D$, paired with some $F'\in F_0\backslash D_{S'}$) and $G\notin D_S$. The product of all $m_F$ will do. Thus $\sigma_S\cap\sigma_{S'}$ is their common face, condition (3) is satisfied.
\end{proof}

\section{More on Sch\"on Varieties}

Our main object here is to prove Theorem \ref{main2}. Let $Y\subset T$ be sch\"on, and $\Delta'$ any fan supported on $\trop(Y)$. Let $\Delta$ be a refinement of $\Delta'$ such that $\Delta$ is sch\"on. Denote $\oY$ and $\oY'$ to be the closure of $Y$ in $X(\Delta)$ and $X(\Delta')$ respectively. We have a commutative diagram,
$$\begin{xymatrix}
{ \oY \ar[r] \ar[d] & X(\Delta) \ar[d] & T\times_k \oY \ar[l] \ar[d] \\
\oY' \ar[r] & X(\Delta') & T\times_k \oY' \ar[l] }
\end{xymatrix}$$
The first step is to show that the left square is a set-theoretic cartesian diagram, and so is the right square (theorem \ref{setinverse}).  

\begin{Lemma}\label{lem1}
Let $Z\subset T$ be a Zariski closed subset of equidimension. If $\trop(Z)$ is preserved by translation by a linear subspace $L$ of $N_{\bbQ}$, then each irreducible component of $Z$ is preserved by the corresponding subtorus of $T$.
\end{Lemma}

\begin{proof}
Let $T\rightarrow T'$ be a homomorphism of tori corresponding to $N_{\bbQ}\rightarrow N_{\bbQ}/L$. Let $Z_i\subset Z$ be irreducible components, and $Y_i$ the closure of the image of $Z_i$ in $T'$. Since $\trop(Y_i)$ is the image of $trop(Z_i)$ in $N_{\bbQ}/L$, we have $\dim\trop(Y_i)\leq\dim\trop(Z)-\dim L=\dim Z-\dim L$. On the other hand $\dim\trop(Y_i)=\dim Y_i\geq \dim Z_i-\dim L$, hence $\dim Z_i=\dim Y_i+\dim L$. Thus the generic fiber of $Z_i\rightarrow Y_i$ is the relative torus, hence $Z_i\cong Y_i\times\ker(T\rightarrow T')$.
\end{proof}

\begin{Lemma}\label{lem2}
Let $\oY\subset X(\Delta)$ be a tropical compatification, $W\subset X(\Delta)$ an orbit closure, and $Z=\oY\cap W$, then $Z$ is equidimensional and each irreducible component intersects the open orbit in $W$.
\end{Lemma}
\begin{proof}
By pulling back the structure map, we see that the structure map of $Z$ in $W$ is also flat. By the openness of flat morphism, each irreducible component of $Z$ intersects the open orbit $O$ of $W$. Note that $(Z\cap O)\times\ker(T\rightarrow O)$ is the fiber of $T\times\oY\rightarrow X(\Delta)$ at any point in $O$, thus $Z\cap O$ is equidimensional and so is $Z$.
\end{proof}

\begin{Lemma}\label{lem3}
Notations and assumptions as in Lemma \ref{lem2}, if $p:X(\Delta)\rightarrow X(\Delta')$ is a proper toric map, then $Z$ is preserved by $T_W$ where $T_W$ is the relative torus of $W\rightarrow p(W)$.
\end{Lemma}

\begin{proof}
Since $Z\subset W$ is a tropical compactification of $Z\cap O$, $\trop(Z\cap O)$ is the support of the fan of $W$ (see remark \ref{importantremark}), which is the inverse image of the support of the fan of $p(W)$ since $p$ is proper. Thus $\trop(Z)$ is preserved by translation of a linear subspace of $N_O$ whose corresponding subtorus of $O$ is the relative torus $T_W$. It then follows from Lemma \ref{lem1}.
\end{proof}

\begin{Theorem}\label{setinverse}
Let $Y\subset T$ be a closed subvariety, $\oY\subset X(\Delta)$ a tropical compactification, $p:X(\Delta)\rightarrow X(\Delta')$ a proper toric map,we have a commutative diagram
$$
\begin{xymatrix}
{ \oY \ar[r] \ar[d] & X(\Delta) \ar[d]^p \\
Y' \ar[r] & X(\Delta'), }
\end{xymatrix}
$$
where $Y'=p(\oY)$, then $\oY$ is a set-theoretic inverse image of $Y'$ under the map $X(\Delta)\rightarrow X(\Delta')$. This applies in particular to the case when $\Delta$ is a refinement of $\Delta'$.
\end{Theorem}
\begin{proof}
We show that if $y\in p^{-1}(y')\cap\oY$ then $p^{-1}(y')\subset\oY$. Let $O'$ be the torus orbit containing $y'$, $p^{-1}(O')$ is a union of toric varieties with connected fibres. It suffices to show that if $y'$ falls in one irreducible component of $p^{-1}(O')$, then the whole fiber in that irreducible component is contained in $\oY$, then this follows from Lemma \ref{lem3}.
\end{proof}

\begin{Corollary}
Notations as in Theorem \ref{setinverse}, let $P$ be the fiber product $Y'\times_{X(\Delta')} X(\Delta)$, then the induced map $\oY\rightarrow P$ is the reduction of $P$.
\end{Corollary}

\begin{proof}
$P$ is a closed subscheme of $X$ and so is $Y$, the induced map $\oY\rightarrow P$ is surjective by Theorem \ref{setinverse}, since $\oY$ is integral, we have $\oY=P_{\red}$.
\end{proof}

The following lemma completes the proof of Theorem \ref{main2}.

\begin{Lemma} \label{lem4}
In the following diagram, $X,X',Y$ and $Y'$ are varieties with $X,X'$ and $Y$ normal, $p$ is proper with connected fibers.  
$$\begin{xymatrix}
{ Y \ar[r]^f \ar[d] & X \ar[d]_{p} \\
Y' \ar[r]^{f'} & X'.}
\end{xymatrix}$$
Let $P=Y\times_{X'}X$, assume the induced map $Y\rightarrow P$ is the reduction of $P$. If $f$ is smooth, then so is $f'$.
\end{Lemma}

\begin{proof}
Let $y'\in Y'(k)$ be a closed point, we show that $f'$ is smooth at $y'$. Let $x'=f'(y')$, $F=Y'_{x'}$, the scheme-theoretic fibre over $x'$ viewed as a closed subscheme of $Y'$, and $G=X_{x'}$. Let $\tY$ be the normalization of $Y$, then the map $Y\rightarrow Y'$ factors through $\tY$, as we have a following diagram:
$$\begin{xymatrix}
{Y \ar[r]^i \ar[d] & P \ar[r]^j \ar[d]_{q} & X \ar[d]_{p} \\
\tY \ar[r] & Y' \ar[r]^{f'} & X'.}
\end{xymatrix}$$

Note that $j^{-1}(G)=F\times_k G$, and $f^{-1}(G)$ is a closed subscheme of $j^{-1}(G)$ with same support, and is smooth over $G$, thus we have $f^{-1}(G)=F_{\red}\times_k G$ and $F_{\red}$ is regular.

Let $f_1,\ldots,f_d\in\cO_{F_{\red},y'}$ be a regular system of parameters, and lift them to the local ring $\cO_{Y',y'}$ with the same notation. Assume $f_1,\ldots,f_d\in\Gamma(U,\cO_{Y'})$, and we may shrink $U$ a little to assume $U\cap V(f_1,\ldots,f_d)\cap F=\{y'\}$.

Let $U\rightarrow X'\times_k\bbA^d_k$ be the map defined by $(f',f_1,\ldots,f_d)$, and $\varphi:Y_U\rightarrow X\times_k\bbA^d_k$ defined by $(f,f_1,\ldots,f_d)$. For any point $y\in Y$ lying over $y'$, $f_1,\ldots,f_d$ restricted in $\cO_{f^{-1}(f(y)),y}$ is a regular system of parameters since $f^{-1}(f(y))=F_{\red}$. By Lemma \ref{lem5}, $\varphi$ is \'etale at $y$, thus $\varphi$ is \'etale in a neighbourhood of $Y_y'$.

Since $p$ has connected fibres, $\tY\rightarrow Y'$ is a homeomorphism of the underlying topological spaces, thus there is a unique $\tilde{y}\in\tY$ lying over $y'$ and $\varphi$ is also \'etale in a neighbourhood of $Y_{\tilde{y}}$. We have $\varphi^{-1}(X_x,0)=Y_y$, so apply Lemma \ref{lem6} for $W=Y_y\subset Y$ and $Z=(X_x,0)\subset X\times_k\bbA^d_k$, and formal function theorem for proper maps $Y\rightarrow\tY$ and $X\times\bbA^d\rightarrow X'\times\bbA^d$, we have isomorphism of formal local rings $\widehat{\cO}_{\tY,\tilde{y}}\cong\widehat{\cO}_{X'\times\bbA^d,(x,0)}$. 

Thus $\tY\rightarrow X'\times\bbA^d$ is \'etale at $\tilde{y}$, so it separates tangent vectors at $\tilde{y}$. It follows that $\tY\rightarrow Y'$ also separates tangent vectors at $\tilde{y}$, hence $\tY\rightarrow Y$ is a closed immersion, we have $\tY\cong Y'$, and $f':Y'\rightarrow X'$ is smooth.

\end{proof}

\begin{Lemma} \label{lem5}
Let $f:Y\rightarrow X$ be a morphism of $k$-schemes, $y\in Y(k)$ a closed point, then $f$ is smooth at $y$ iff there exist $f_1,\ldots,f_d\in\fm_y\subset\cO_{Y,y}$ such that the locally defined map $Y\rightarrow X\times_k \bbA^d_k$ given by $(f,f_1,\ldots,f_d)$ is \'etale at $y$.
\end{Lemma}

\begin{proof}
If there is a map $Y\rightarrow X\times_k\bbA^d_k$ is \'etale at $y$, by composing a smooth map $X\times_k\bbA^d_k\rightarrow X$, $Y\rightarrow X$ is smooth at $y$. Conversely we know that $Y_x$ is regular at $y$, choose $f_1,\ldots,f_d$ in $\cO_{Y,y}$ which form a regular system of parameters in $\cO_{Y_x,y}$, consider the map
$$\varphi:\cO_{X,x}[X_1,\ldots,X_d]_{(X_1,\ldots,X_d)}\rightarrow\cO_{Y,y}$$
given by $X_i\mapsto f_i$. $\varphi$ is flat since $\cO_{X,x}\rightarrow \cO_{X,x}[X_1,\ldots,X_d]_{(X_1,\ldots,X_d)}$ is faithfully flat, and $\varphi$ is also geometrically regular since the geometric fiber of the closed point is a reduced point, thus $\varphi$ is smooth of relative dimension 0, hence \'etale.
\end{proof}

\begin{Lemma}\label{lem6}
Let $f:Y\rightarrow X$ be an \'etale morphism of schemes, $Z\subset X$ a closed subscheme, if $W=f^{-1}Z\rightarrow Z$ is an isomorphism, then $W_n\rightarrow Z_n$ is an isomorphism for all $n>0$. 
\end{Lemma}

\begin{proof}
Let $I$ and $J$ be the ideal sheaf of $Z$ and $W$ respectively. Clearly $J=f^{-1}I\cdot\cO_{Y}$, and $J^n=f^{-1}I^n\cdot\cO_Y$, thus $W_n=f^{-1}(Z_n)$ and we have a morphism $W_n\rightarrow Z_n$ for each $n>0$ and it is an isomorphism for $n=1$. It follows that $W_n\rightarrow Z_n$ is a homeomorphism for the underlying spaces, and is \'etale. $Z$ is defined by a nilpotent ideal in $Z_n$, by the formal property of \'etale, 
$$\begin{xymatrix}
{Z=W \ar[r] \ar[d] & W_n \ar[d] \\
Z_n \ar[r]^{=} \ar@{-->}[ru]& Z_n,}
\end{xymatrix}$$
$W_n\rightarrow Z_n$ admits a section, then it follows that this is an isomorphism.

\end{proof}

\section{Toric Schemes over a Discrete Valuation Ring}
In order to extend tropical compactification to non-constant coefficient case, we need to introduce the notion of toric schemes over a discrete valuation ring. This has first appeared in Mumford's \cite{KKMSD} in an effort of extending the semi-stable reduction theorem. Smirnov \cite{Smi97} constructed toric schemes from certian polyhedra (his so-called admissible polyhedra). We observe that a simple construction is just by pulling back a toric map $X(\Delta)\rightarrow \bbA^1_k$ to the spectrum of the local ring at 0. We shall also compare a toric scheme (with same combinatorial data) in various ring extensions since we ultimately intend to study a subvariety $Y\subset T$ over the field of Puiseux series.

First let's fix some notations, $(R,\fm,k)$ is a discrete valuation ring with $\fm$ the maximal ideal, $k$ the residue field and $K$ the quotient field, assume $k\subset R\subset K$. Let $\bbK$ be the algebraic closure of $K$. We write $\Spec R=\{\eta,s\}$ where $\eta$ is the generic point and $s$ the closed point. For any scheme $X$ over $\Spec R$, we use $X_{\eta}$ and $X_s$ to denote the generic fibre and special fiber respectively. We fix $t\in\fm\backslash\fm^2$, a uniformizer. Our main application would be that $k$ is an algebraically closed field of characteristic 0, $\bbK$ is the algebraic closure of $k(t)$, i.e. the field of Puiseux series over $k$, $K=k((t^{1/n}))$ for some $n>0$.

\begin{Definition}\label{defofTscheme}
A (normal) toric scheme over $R$ is an integral normal scheme $\fX$ together with a map $\fX\rightarrow\Spec R$ which is separated and of finite type, such that it contains the torus $T_K$ in its generic fibre $\fX_{\eta}$ and there is a group scheme action $T_R\times_R\fX\rightarrow\fX$ which extends the left multiplication of $T_K$ on itself.
\end{Definition}

To construct a toric scheme, we need the following data. Let $N$ be a lattice of rank $n$, $\tN=N\oplus\bbZ$, $\Delta\subset\tN_{\bbQ}$ a fan such that $\pr_2(\Delta)\subseteq \bbQ_{\geq 0}$. We call such a fan an \emph{admissible fan} and cones in it \emph{admissible cones}. This gives a toric morphism of ordinary toric varieties $X(\Delta)\rightarrow\bbA^1_k$, the toric scheme $\fX(\Delta)\rightarrow\Spec R$ is obtained from $X(\Delta)$ via the base change $\Spec R\rightarrow\bbA^1_k$, where the map $\Spec R\rightarrow\bbA^1_k$ is defined by the map $k[x]\rightarrow R:x\mapsto t$. If $\Delta$ consists of a single cone $\sigma$ together with all its faces, we write $\fX_{\sigma}$ instead of $\fX(\Delta)$.

This construction is essentially the same as in \cite{KKMSD}. We shall describe $\cO(\fX_{\sigma})$ explicitly for $\sigma\in\Delta$. Let $M$ and $\tM =M\oplus\bbZ$ be the dual of $N$ and $\tN$ respectively. For a lattice point $\tm\in\tM$, we write $\tm=(m,s)$ where $m\in M$ and $s\in\bbZ$, let $e=(0,1)\in\tM$. Note that for any admissible cone $\sigma$, $e$ is contained in its dual cone $\sigma^{\vee}$. If we write $\Ssig$ for the monoid $\sigma^{\vee}\cap\tM$, then
$$\cO(\fX_{\sigma})=k[\Ssig]\otimes_{k[x]}R=R[\chi^{\tm}]_{\tm\in\Ssig}/(\chi^e-t)=R[\chi^m t^s]_{(m,s)\in\Ssig}.$$
Denote $R[\chi^m t^s]_{(m,s)\in\Ssig}$ by $A[\Ssig]$, note that $A[S]$ can be defined for any submonoid $S\subset \tM$ containing $e$. In \cite{KKMSD}, $\fX(\Delta)$ is constructed by glueing all $\Spec A[\Ssig]$, we refer the reader to \cite{KKMSD} IV \S3 for basic properties of $\fX(\Delta)$ which are analogous to ordinary toric varieties, for example one-to-one correspondence of toric orbit and cones, correspondence of $T$-linearized invertible sheaves and piecewise linear maps $|\Delta|\rightarrow \bbQ$, etc. We give a brief proof of the following properties of $\fX(\Delta)$ since they are not included in \cite{KKMSD} and are sort of special for toric $R$-schemes.

\begin{Proposition}
\begin{enumerate}
\item The generic fiber of $\fX$ is $X(\Delta\cap N)$, where $\Delta\cap N$ is the subfan of $\Delta$ consisted of all cones of $\Delta$ which are in $N\subset\tN$.
\item For any $\rho\in\Delta(1)$, a 1-dimensional cone, let $v_{\rho}$ be the first lattice point on $\rho$, then $\val_{V_{\rho}}(\chi^m t^s)=\left\langle \tm,v_{\rho}\right\rangle$, where $\tm=(m,s)\in\tM$ and $V_{\rho}$ is the closure of the orbit $O_{\rho}$ corresponding to $\rho$.
\item The special fiber of $\fX(\Delta)$ is a union of irreducible components indexed by $\rho\in\Delta(1)$ not contained in $N$, whose reduction is the union of $V_{\rho}$. The irreducible component indexed by $\rho$ is reduced if and only if $\pr_2(v_{\rho})=1$, the special fiber $\fX_s$ is reduced iff $\pr_2(v_{\rho})=1$ for all $\rho\not \subset N$.
\item For an admissible fan $\Delta\subset\tN=N\oplus\bbZ$, let $\Delta[d]$ be the same fan but considered in $N\oplus d\bbZ$. There exists an integer $d>0$ such that for any discrete valuation ring $R'$, $\fX(\Delta[d])$ as a toric $R'$-scheme has reduced special fiber, and furthermore when $R'=R[t^{1/d}]$, $\fX(\Delta[d])$ is obtained from $\fX(\Delta)$ via base change $\Spec R[t^{1/d}]\rightarrow\Spec R$ followed by normalization.
\end{enumerate}
\end{Proposition}

\begin{proof}
(1) follows easily from the pullback construction of $\fX$ and the toric orbit structure of the ordinary toric variety $X(\Delta)$.

Let $\Osig'$ and $\Vsig'$ be the usual toric orbit and orbit closure of $X(\Delta)$, it is easy to check that $\Osig$ and $\Vsig$ are obtained via base change $\fX\rightarrow X(\Delta)$. For (2), note that the rational function $\chi^m t^s\in K[M]$ on $\fX$ is induced from the rational function $\chi^{\tm}\in k[\tM]$ on $X(\Delta)$, thus $\val_{V_{\rho}}\chi^m t^s= \val_{V_{\rho}}\chi^{\tm}=\left\langle\tm,v_{\rho}\right\rangle$.

For (3), we already know $(\fX_s)_{\red}=\bigcup_{\rho} V_{\rho}$ for $\rho\not\subset N$. Let $\fX_{s,\rho}$ be the irreducible component indexed by $\rho$, with $V_{\rho}$ its reduction. $\fX_{s,\rho}$ is reduced iff $\val_{\Vsig} t=1$, that is $\pr_2(v_{\rho})=1$ by (2). This can also be shown by the affine ring, let $\sigma\in\Delta$, then $\cO(\fX_{\sigma}\times_R k)=k[\sigma^{\vee}\cap \tM]/(\chi^e)$, the effect of modulo $(\chi^e)$ is to make $\chi^{\tm}=0$ if $\tm-e\in\sigma^{\vee}\cap\tM$. If $\pr_2(v_{\rho})=0$ or 1, for all $\rho<\sigma$, then $\chi^{\tm}=0$ for all $\tm$ in the interior of $\sigma^{\vee}$, hence $k[\sigma^{\vee}\cap \tM]/(\chi^e)$ is reduced.

For (4), let $d$ be the least common multiple of $v_{\rho}$ for all $\rho\not\subset N$, then $\Delta[d]$ satisfies the assumption in (3), hence $\fX(\Delta[d])$ has reduced special fiber. Now assume $R'=R[t^{1/d}]$, we show that $\fX(\Delta[d])$ is the normalization of $\fX(\Delta)\times_R R[t^{1/d}]$. This is a local question, we show this for each $\sigma\in\Delta$. Note that the dual of $N\oplus d\bbZ$ is $M\oplus\frac 1 d \bbZ$ and we have a natural embedding $\tM\subset M\oplus\frac 1 d\bbZ$ coresponding to $N\oplus d\bbZ\subset\tN$. Let $\sigma[d]$ be the same cone but considered in $N\oplus d\bbZ$, and $S_{\sigma[d]}$ the monoid $\sigma^{\vee}\cap(M\oplus\frac 1 d\bbZ)$.

The affine ring of $\fX_{\sigma}\times_R R[t^{1/d}]$ is $R[\chi^m t^s]_{(m,s)\in\Ssig}\otimes_R R[t^{1/d}]$, the effect of this tensor product is just introducing $\frac 1 d e$ in $\Ssig$ i.e. if $S'$ is the monoid generated by $\Ssig$ and $\frac 1 d e$ in $S_{\sigma[d]}$, then 
$$R[\chi^m t^s]_{(m,s)\in\Ssig}\otimes_R R[t^{1/d}]=R[\chi^m t^s]_{(m,s)\in S'}.$$
The effect of normalization is saturating $S'$, which is $S_{\sigma[d]}$.
\end{proof}

\begin{Remark}\label{rmkfx}
If $\fX(\Delta)$ has reduced special fibre, then $\fX(\Delta[d])=\fX(\Delta)\times_R R[t^{1/d}]$, this is because $S'$ (the monoid generated by $\Ssig$ and $\frac 1 d e$ in $S_{\sigma[d]}$) is already saturated, i.e. $S'=S_{\sigma[d]}$.
\end{Remark}

\section{Extension to Non-constant Coefficient Case} 
We first define tropical and sch\"on compactification for $Y\subset T$ over a field with discrete valuation and establish parallel results as in the constant coefficient case. Notations as in the previous section, let $Y\subset T_K$ be a subvariety, $\trop(Y)$ in this case is taken to be the tropicalization of $Y_{\bbK}\subset T_{\bbK}$.

Note that $\trop(Y)$ sits inside $N_{\bbQ}$, define $\cT(Y)$ to be the closure of the set $\{(t\cdot x,t)\in N_{\bbQ}\oplus\bbQ=\tN_{\bbQ}|x\in\trop(Y),t\in\bbQ_{>0}\}$. We take $\cT(Y)$ as a replacement of $\trop(Y)$ in the non-constant coefficient case. If we modify the BG-set definition of $\trop(Y)$ (definition \ref{tropset}, (2)), $\cT(Y)$ coincides with the following set:
$$\{(u|_{M_{\bbQ}}, u(t))\in N_{\bbQ}\oplus\bbQ|u:\cO(Y)\rightarrow \bbQ\}$$ 
where $u$ runs over all valuations trivial on $R^*$ and non-negative on $R$. Let $\Delta\subset\tN_{\bbQ}$ be an admissible fan, $\fX(\Delta)$ the corresponding toric scheme over $R$, $\oY$ the closure of $Y$ in $\fX(\Delta)$, we make the following definitions.

\begin{Definition}
We say $\oY$ is a \emph{tropical compactification} or $\Delta$ is a \emph{tropical fan} if $\oY$ is proper over $R$ and the structure map $T_R\times_R\oY\rightarrow \fX(\Delta)$ is flat and surjective.
\end{Definition}

\begin{Definition}
$\oY$ is called a \emph{sch\"on compactification} if it is tropical and moreover the structure map is smooth. We say $Y$ is \emph{sch\"on} in $T$ if it admits a sch\"on compactification, and we say $Y$ is \emph{sch\"on} if it is sch\"on in the intrinsic torus.
\end{Definition}

\begin{Definition}
$\oY$ is called a \emph{h\"ubsch compactification} if it is sch\"on, $\fX(\Delta)$ has reduced special fibre and $K_{\oY}+B_{\oY}$ is ample. We say $Y$ is \emph{h\"ubsch} in $T$ if it admits a h\"ubsch compactification, and we say $Y$ is \emph{h\"ubsch} if it is h\"ubsch in the intrinsic torus.
\end{Definition}

\begin{Remark}
Since $\fX(\Delta)$ has reduced fiber and $\oY$ has smooth structure map, $K_{\fX}$ and $K_{\oY}$ are well defined. We have $K_{\fX}+B_{\fX}=0$, and by adjunction formula, $K_{\oY}+B_{\oY}=\det\sN_{\oY/\fX(\Delta)}$.
\end{Remark}

\begin{Proposition}\label{supportontrop}
$\oY$ is proper over $R$ iff $|\Delta|\supset\cT(Y)$. If $\Delta$ is a tropical fan, then $|\Delta|=\cT(Y)$.
\end{Proposition}

\begin{proof}
This can be proved following the same idea as in the constant coefficient case in \cite{Tev07}, here we proceed with a new proof using the BG-set definition and valuative criterion.

Suppose $\oY$ is proper over $R$, let $u$ be a valuation $K(Z)^{\times}\rightarrow\bbQ$ which is trivial on $R^*$ and nonnegative on $R$ for some subvariety $Z$ of $Y$. Let $R(Z)$ be the valuation ring, then we have a following commutative diagram,
$$\xymatrix{\Spec K(Z) \ar[r]\ar[dd]  \ar[r] & \oY \ar[r] & \fX(\Delta) \ar[dd] \\
&  \fX_{\sigma} \ar[ur] \\
\Spec R(Z)\ar[rr] \ar[uur]^f \ar[ur]_g & & \Spec R\,,
}$$
where the existence of $f$ is due to the valuative criterion of properness of $\oY$, $f$ then factors through some open affine toric scheme $\fX_{\sigma}\subset\fX(\Delta)$. $g$ corresponds to a ring homomorphism $\cO(\fX_{\sigma})\rightarrow R(Z)$, which implies that $u$ takes non-negative values on $\Ssig$, thus $[u]\in\sigma$ where $[u]$ denotes the image of $u$ in $\tN_{\bbQ}$. We proved $|\Delta|\supset \cT(Y)$.

Suppose $|\Delta|\supset \cT(Y)$, if $\oY$ is not proper, there is a proper $R$-scheme $Y'$, containing $\oY$ as an open dense subscheme. Let $y_1$ be the generic point of $Y'$ and $y_0$ any point of $Y'$ not contained in $\oY$, then $\cO_{Y',y_0}$ is dominated by a discrete valuation ring of $K(Y')=K(Y)$. Let $u$ denote the valuation, and $R(Y)$ the valuation ring, then $[u]\in\sigma$ for some $\sigma\in\Delta$. The map $\Spec R(Y)\rightarrow Y'$ sending the generic point to $y_1$ and the closed point to $y_0$ factors through $\fX_{\sigma}$, a contradiction.

Let $\Delta$ be a tropical fan, we already proved $\Delta\supset \cT(Y)$, hence for any $\Delta'$ refining $\Delta$, $\oY(\Delta')$ is proper over $R$.
Let $\fY=\oY(\Delta)\times_{\fX(\Delta)}\fX(\Delta')$. We show that $\fY=\oY(\Delta')$. Indeed we have a fibre diagram
$$\xymatrix{T_R\times_R \fY \ar[r]^f \ar[d] \ar@{}|{\Box}[dr] & \fX(\Delta ') \ar[d] \\
T_R\times_R \oY(\Delta) \ar[r] & \fX(\Delta)\,. \\}$$
Since $\oY(\Delta)$ is tropical, the top arrow is also flat and surjective. Restricting on the open subscheme $T_K\subset\fX(\Delta')$, $f^{-1}(T_K)=T_K\times_K Y$ which is integral. By the lemma below, $T_R\times_R\fY$ is integral, so is $\fY$. Hence $\fY$ is $\pi^{-1}(\oY(\Delta))\subset\fX(\Delta')$ with reduced induced structure, hence $\oY(\Delta')=\fY$, it is tropical and is the pullback of $\oY(\Delta)$.

It remains to prove that $\Delta$ is supported on $\cT(Y)$. Suppose it is not, let $\Delta'$ refines $\Delta$ such that there is a subfan $\Delta''\subset\Delta'$ with $|\Delta''|=\cT(Y)$. $\fX(\Delta'')$ is a open subscheme of $\fX(\Delta')$, but not equal to $\fX(\Delta')$. $\oY(\Delta'')=\oY(\Delta')$, thus the structure map $T_R\times_R \oY(\Delta')\rightarrow\fX(\Delta')$ fails to be surjective. 
\end{proof}

\begin{Lemma}\label{integrallemma}
Let $f:X\rightarrow Y$ be a flat morphism of schemes. Assume $Y$ is integral. If there exists a dense Zariski open subset $U\subset Y$ such that $f^{-1}(U)$ is integral, then $X$ is integral. 
\end{Lemma}

\begin{proof}
For any $x\in X$, let $y=f(x)$. Take affine open neighbourhoods $\Spec B$ and $\Spec A$ of $x$ and $y$ respectively such that $f:\Spec B\rightarrow\Spec A$, then $A\rightarrow B$ is flat. We can find $a\in A$ such that $\Spec A_a\subset U$. Since $f^{-1}(\Spec A_a)\cap \Spec B=\Spec B_a$, we know $B_a$ is a domain. Tensoring $0\rightarrow A\rightarrow A_a$ with $B$ over $A$, we have $0\rightarrow B\rightarrow B_a$, hence $B$ is a domain. $X$ is reduced.

If $X$ is not irreducible, there is an irreducible component $X'$ of $X$ such that $f:X'\rightarrow Y-U$. However each irreducible component of $X$ should dominate $Y$ by the openness of a flat map, this is a contradiction, which proves the lemma.
\end{proof}

\begin{Proposition}
Any refinement of a tropical fan is tropical, if $Y\subset T_K$ admits a sch\"on compactification, then any tropical fan produces a sch\"on compactification. Let $\oY\subset\fX(\Delta)=\fX$ be a sch\"on compactification, then $\oY$ is locally a complete intersection in $\fX(\Delta)$. Furthermore for any refinement $\Delta'$ of $\Delta$, let $\oY'$ be the closure of $Y$ in $\fX(\Delta')=\fX'$ with proper birational map $\pi\,:\oY'\rightarrow\oY$, then $\det\sN_{\oY'/\fX'}=\pi^*(\det\sN_{\oY/\fX})$.
\end{Proposition}
\begin{proof}
Same proof as in the constant coefficient case \cite{Tev07} using the fibre diagram in proposition \ref{supportontrop}.
\end{proof}

\begin{Proposition}
Tropical fan exists assuming $Y\subset T$ is rigid.
\end{Proposition}

\begin{proof}
Same proof as in the constant coefficient case using Kapranov's visible contour construction and Laforgue's transversality argument.
\end{proof}

\begin{Lemma}\label{preserveunderextension}
Let $Y\subset T_K$ be a subvariety, assume $Y$ is geometrically integral, if $\Delta$ produces a tropical (resp. sch\"on, resp. h\"ubsch) compactification of $Y$, then $\Delta[d]$ also produces a tropical (resp. sch\"on, resp. h\"ubsch) compactification for $Y_{K'}\subset T_{K'}$ where $K'=K(t^{1/d})$.
\end{Lemma}

\begin{proof}
Let $R'=R[t^{1/d}]$, $\oY(\Delta[d])$ be the closure of $Y_{K'}$ in $\fX(\Delta[d])$, we have a fibre diagram (by lemma \ref{integrallemma})
$$\xymatrix {T_{R'}\times_{R'}\oY(\Delta[d]) \ar[r] \ar[d] & \fX(\Delta[d]) \ar[d] \\ T_R\times_R \oY \ar[r] & \fX(\Delta).}$$
Hence if $\oY$ is tropical (resp. sch\"on), so is $\oY(\Delta[d])$.

If $\oY$ is h\"ubsch, $\fX(\Delta)$ has reduced special fibre and so does $\fX(\Delta[d])$, thus $\oY(\Delta[d])$ and $\fX(\Delta[d])$ are obtained from $\oY$ and $\fX(\Delta)$ simply by ring extension $\Spec R'\rightarrow \Spec R$ (see remark \ref{rmkfx}), therefore $\oY(\Delta[d])$ is also h\"ubsch.

\end{proof}

\begin{Proposition}
If $Y$ is h\"ubsch in $T_K$, assume $Y$ is geometrically integral, then $\cT(Y)$ has a minimal fan stucture $\Delta$ corresponding to the log canonical compactification.
\end{Proposition}

\begin{proof}
Let $\Delta'$ be another fan supported on $\cT(Y)$, we show that $\Delta'$ is a refinement of $\Delta$. Suppose on the contrary $\Delta'$ does not refine $\Delta$, then there exists a $d$-dimensional cone $\sigma'\in\Delta'$ ($d=\dim Y$), which is not contained in any cone of $\Delta$. There is a $(d-1)$-dimensional cone $\alpha\in\Delta$ which meets the interior of $\sigma'$. Let $\Delta''$ be a common strictly simplicial refinement, there is a $(d-1)$-dimensional cone $\alpha''\subset\alpha\cap\sigma'$, meeting the interior of both $\alpha$ and $\sigma'$. We may assume $\fX(\Delta),\fX(\Delta')$ and $\fX(\Delta'')$ all have reduced fibre, otherwise consider $\Delta[l],\Delta'[l]$ and $\Delta''[l]$ for some $l$ (lemma \ref{preserveunderextension}).

We have proper birational maps 
$$p_1:\,\fX'':=\fX(\Delta'')\rightarrow\fX:=\fX(\Delta)$$
and 
$$p_2:\,\fX''\rightarrow\fX':=\fX(\Delta').$$
Let $Z$, $Z'$ and $Z''$ be the orbit closure $V_{\alpha}$, $V_{\sigma'}$ and $V_{\alpha''}$ in $\fX$, $\fX'$ and $\fX''$ respectively.

Note that $\sigma'$ is not contained in $N$ since it's a maximal cone in $\cT(Y)$, and consequently same for $\alpha$ and $\alpha''$, hence $Z$, $Z'$ and $Z''$ are all on the special fiber, in particular, are all normal toric varieties over $k$.

$Z'$ is isomorphic to a torus $T'_k$. The induced morphism $p_2:\,Z''\rightarrow Z'$ is proper toric morphism of relative dimension 1, thus $Z''\cong Z'\times_k\bbP^1_k$. $p_1:\, Z''\rightarrow Z$ is birational.

Let $\oY$ and $\oY''$ be the closure of $Y$ in $\fX$ and $\fX''$ respectively. $\oY''$ is a sch\"on compactification with reduced special fibre. The scheme-theoretic intersection $C:=Y''\cap Z''\subset T'_k\times_k\bbP^1_k$ is 1-dimensional, reduced, proper and smooth over $k$, thus $C\cong z\times\bbP^1$ for some $z\in T'_k$ a 0-dimensional reduced closed subscheme. So we have that $K_C+B_C$ is trivial. Since $\oY''\rightarrow\oY$ is log crepant, and $K_{\oY}+B_{\oY}$ is ample, by projection formula we conclude that $C$ is contracted by $p_1:\,Z''\rightarrow Z$. Since $p_1$ is equivariant, all fibers $z'\times\bbP^1$ are contracted. A contradiction against that $p_1$ is a birational map.
\end{proof}

\begin{Theorem}
Assume $\charact K=0$, then for any variety $Y$ over $K$, $Y$ contains a sch\"on very affine variety.
\end{Theorem}

\begin{proof}
The proof proceeds in an analogous streamline as in the constant coefficient case , only with some technical issues to be taken care of. Let $Y^o$ be a regular very affine variety over $K$ and $\oY$ a regular compactification, projective over $\Spec R$ with reduced special fiber and s.n.c. boundary divisor. Let $M$ be a lattice with a homomorphism $\phi:M\rightarrow \cO^*(Y)$ such that $K[M]\rightarrow \cO(Y)$ is surjective. For any stratum $S$, let $\tM_S$ be the set $\{\tm=(m,r)\in\tM=M\oplus\bbZ|\tm\in\cO(\Star(S))\}$ (where we think $\tm=\phi(m)t^r$ as a rational function on $Y$), the conditions in proposition \ref{schoncondition} are now the following: 
\begin{enumerate}

\item for each strata $S$, $\Star M$ is affine and $A[M_S]\rightarrow \cO(Star(S))$ is surjective.
\item for each strata $S$ and any divisor $D_0\in D_S$, there exists $\tm\in \tM_S$ such that $\val_{D_0}\tm=1$ and $\val_{D'}\tm=0$ for all $D'\in D_S\backslash\{D_0\}$. 
\item the cones $\sigma_S$ (the dual cone of $\tM_S$) form an admissible fan in $\tN_{\bbQ}$.
\end{enumerate}

We may assume $Y$ is regular, by Hironaka's resolution theorem and Mumford's semistable reduction theorem, there is a compactification $\oY$, possibly over a ring extension $R\subset R[t^{1/d}]$, which is regular with reduced special fibre and s.n.c. boundary divisor. Thanks to the following lemma of a relative version of Bertini's theorem, the above conditions can be achieved by adding more generic hyperplane sections as in the proof of theorem \ref{main1}. 
\end{proof}

\begin{Lemma}[\cite{JS07}]
$X$ is regular scheme, flat and quasi-projective over $\Spec R$, assume $X_s$ is reduced and s.n.c., then a general hyperplane $H\subset \bbP_R^n$ intersects $X$ transversely and $(X\cap H)\cup X_s$ is s.n.c.
\end{Lemma}

\begin{Theorem}
If $Y$ is sch\"on in $T_K$, then any fan supported on $\cT(Y)$ produces a sch\"on compactification.
\end{Theorem}

\begin{proof}
This is a relative version of theorem \ref{main2}, which is proved by reducing to the constant case. Let $\Delta'$ be any fan supported on $\cT(Y)$, and $\Delta$ is a refinement of $\Delta'$ which is sch\"on, let $\oY$ and $\oY'$ be the closure of $Y$ in $\fX(\Delta)$ and $\fX(\Delta')$ respectively. We have a commutative diagram as follows:
$$\xymatrix{ T_R\times_R\oY \ar[r]^f  \ar[d] & \fX(\Delta) \ar[d]^p \\
T_R\times_R\oY' \ar[r]^g & \fX(\Delta'). } $$

Note that $f$ is smooth. Pull the diagram back to the generic fibre and lift to the algebraic closure $\bbK$, we have
$$\xymatrix{ T_{\bbK}\times_{\bbK} \oY_{\bbK} \ar[r]^{f_{\bbK}} \ar[d] & \fX(\Delta)_{\bbK} \ar[d] \\
T_{\bbK}\times_{\bbK} \oY'_{\bbK} \ar[r]^{g_{\bbK}} & \fX(\Delta')_{\bbK}. } $$
Now $\oY_{\bbK}$ is a disjoint union of irreducible components of pure dimension, each irreducible component has smooth surjective structure map possibly in a smaller toric open set (see remark \ref{importantremark}). Apply theorem \ref{setinverse} to each irreducible component of $\oY_{\bbK}$, we see that $\oY_{\bbK}$ is a set theoretic inverse image of $\oY'_{\bbK}$, hence $\oY'_{\bbK}$ is also a disjoint union of irreducible components of pure dimension. We can now apply theorem \ref{main2} to each irreducible component of $Y_{\bbK}$, $g_{\bbK}$ is smooth.

Let $W$ be any toric orbit closure of $\fX(\Delta)$ on the special fibre (which is a normal toric variety over $k$), let $O\subset W$ be the open orbit, then $W\rightarrow p(W)$ is a proper toric map. Restricting the diagram on $W$, we have 
$$\xymatrix{ \oY\cap W \ar[r] \ar[d] & W \ar[d]^p \\
\oY'\cap p(W) \ar[r] & p(W). }$$
$\oY\cap W$ has smooth structure map in $W$, by the same argument as above, we can show that $\oY'\cap p(W)$ also has smooth structure map. Combining the above result, we see that $f$ is smooth.
\end{proof}

Assume $k$ is of characteristic 0, $\bbK$ is then the field of Puiseux series over $k$. For a subvariety $Y\subset T_{\bbK}$, $Y$ is defined over $K_n:=k((t^{1/n}))$ for some $n$. We propose the following definition.

\begin{Definition}
$Y\subset T_K$ is sch\"on (resp. h\"ubsch) if for some $n$, $Y$ is defined over $K_n$ and $Y\subset T_{K_n}$ is sch\"on (resp. h\"ubsch).
\end{Definition}

This is to avoid some technical difficulty such as schemes over a non-noetherian ring. It is clear that if $Y\subset T_{\bbK}$ is h\"ubsch, $\cT(Y)$ has a minimal fan structure.

\bibliographystyle{amsplain}	
\bibliography{mybib}

\end{document}